\documentclass[preprint,12pt]{elsarticle}

\usepackage{graphicx}

\usepackage{amssymb}
\usepackage{amsthm}
\usepackage{amsmath}
\usepackage{xcolor}

\usepackage[utf8]{inputenc}
\usepackage{hyperref}
\usepackage{natbib}

\journal{Computers \& Fluids}

\begin{document}

\begin{frontmatter}



\title{Brownian dynamics of rigid particles in an incompressible fluctuating fluid by a meshfree method }

\author[AP]{Anamika  Pandey\corref{cor1}}
\ead{pandey@mathematik.uni-kl.de}
\author[SH]{Steffen Hardt}
\ead{hardt@csi.tu-darmstadt.de}
\author[AP]{Axel Klar}
\ead{klar@mathematik.uni-kl.de}
\author[AP]{Sudarshan Tiwari}
\ead{tiwari@mathematik.uni-kl.de}
\address[AP]{Fachbereich Mathematik,  Technische Universit\"{a}t Kaiserslautern, Postfach 3049, 67653  Kaiserslautern, Germany}
\address[SH]{Institute for Nano- and Microfluidics, Center of Smart Interfaces, Technische Universit\"{a}t Darmstadt, Petersenstra$\beta$e 17, 64287 Darmstadt, Germany}
\begin{abstract}
A meshfree Lagrangian  method for the fluctuating hydrodynamic equations (FHEs) with  fluid-structure interactions is presented. Brownian motion of the particle is investigated  by direct numerical simulation of the fluctuating hydrodynamic equations. In this framework a bidirectional coupling has been introduced between the fluctuating fluid and the solid object. The force governing the motion of the solid object is solely due to the surrounding fluid particles. Since a  meshfree formulation is used,  the method  can be extended to  many real applications involving  complex fluid flows. A three-dimensional implementation  is  presented. In particular,  we observe the short and long-time behaviour of the velocity autocorrelation function (VACF) of Brownian particles and compare it with the analytical expression. Moreover, the Stokes-Einstein relation is reproduced to ensure the correct long-time behaviour of Brownian dynamics.  
\end{abstract}

\end{frontmatter}


\section{Introduction}
\label{Int}
The dynamics of small rigid particles immersed in a fluid presents an important and challenging problem, 
in particular,  for  micro/nano scale objects  in small scale geometries. The dynamics of small rigid particle can be influenced by the inherent thermal fluctuation in the fluid. As one approaches  smaller scales, thermal fluctuations play an  essential role in the description of the fluid flow, see for example \cite{oritz2006, Landau1980} or 
\cite{Vailati2011, Donev2011, Francois2012} for more recent works.

This study focuses on  Brownian motion of particles immersed in an incompressible fluid. The average motion of the surrounding fluid yields a hydrodynamic force on the  particles. Moreover, a random force is also experienced by the immersed particles due to the thermal fluctuation in the fluid. The average motion of fluid is modelled by the Navier-Stokes equations. 
The  thermal fluctuations can either be described on a microscopic level using methods like molecular dynamics or they can be included in the continuum description of the fluid by additional stochastic fluxes.  
If one concentrates on a continuum field description, 
the resulting equations of motion for the fluctuating fluid turn out to be  stochastic partial differential equations (SPDEs). Such equations, including an additional stochastic stress tensor in the Navier-Stokes equations, have been proposed by Landau and Lifshitz \cite{Landau1959}.  These equations are termed the 
Landau-Lifshitz Navier-Stokes (LLNS) equations.
Initially, the LLNS equations have been presented for fluctuations around an equilibrium state of the system, but later on, their validity for non-equilibrium systems has also been shown  \cite{espanol1998} and verified by molecular simulations \cite{garcia1991, mansour1987}.

Early work in the context of numerical approximation of the LLNS equations has been done by Garcia et al. \cite{garcia1987}. The authors have developed a simple scheme for the stochastic heat conduction equation and the linearized one-dimensional LLNS equations. Later on in \cite{garcia2007} a centered scheme based on a finite-volume discretization, combined with the third-order Runge-Kutta (RK3) temporal integrator, has been introduced for the compressible LLNS equations. Afterwards, a systematic approach  for the analysis of this grid based finite-volume approximation for the LLNS equations and related SPDEs has been discussed by Donev et al. \cite{garcia2010}. The extension of this numerical solver for the LLNS equations to binary mixtures and staggered schemes for the fluctuating hydrodynamic equations have been presented in \cite{bell2010, bell2012}. A meshfree Lagrangian formulation for the 1D LLNS equations for compressible fluids has been presented by the present authors in \cite{pandey2012} and the results have been compared to the above-mentioned FVM-based RK3 scheme from  \cite{garcia2007}.

In the context of fluid-structure interactions,  Brownian dynamics of immersed particles due to the surrounding fluctuating fluid has been studied. 
  in \cite{patankar2004, sharma2006}, using the coupling of the equations of motion for the  particles with the fluctuating hydrodynamic equations. There,  inertia terms in the governing equations have been neglected and  the resultant time independent problem has been solved with a numerical approach using a fixed grid spatial discretization. A hybrid Eulerian-Lagrangian approach for the inertial coupling of point particle with fluctuating compressible fluids has been presented by \citet{Usabiaga2013}. Subsequently, an inertial coupling method for particles in an incompressible fluctuating fluid has been reported in \cite{Donev2014}. In this work,  the equations of motion of the suspended particle are directly coupled with an incompressible finite-volume solver for the LLNS equations \cite{bell2012}. The authors  have also discussed the Stokes-Einstein relation for fluid-structure systems at moderate Schmidt number, see \cite{Donev2013}. They have modeled the particle through a source term in the momentum equation  and dealt with the full incompressible fluctuating hydrodynamic equations. In this work, a very efficient coupling of the fluctuating fluid with particles and a finite-volume approximation of the coupled model have been proposed. The short and long-time behavior of Brownian dynamics has been captured very well. Moreover, the authors were   able to handle a wide range of Schmidt numbers in their proposed approximation, which has been a difficult task for many numerical approximations. An immersed boundary approach has been reported by \citet{atzberger2011} for fluid-structure interaction with thermal fluctuations using a grid-based method and extended to complex geometries in \cite{atzberger2014}.
The fluctuating hydrodynamics approach has also been used  to analyse the Brownian motion of nanoparticles in an incompressible fluid,   compare \citet{Swaminathan2011}.

The present work distinguishes itself from the existing literature in its approach. An explicit coupling has been used between the fluctuating fluid and the solid structure, and a  numerical approximation based on a meshfree formulation is used for the LLNS equations. In general, meshfree methods are an alternative to classical methods  for problems with time-varying fluid domains such as problems with bodies suspended in a fluid, where one  can avoid  re-meshing during the time evolution.
We note that  a meshfree method termed ``Smoothed dissipative particle dynamics (SDPD)" has been presented in \cite{espanol2003} which incorporates thermal fluctuations.
 The SDPD is a combination of meshfree smoothed particle hydrodynamics (SPH) \cite{liu2003} and dissipative particle dynamics (DPD) \cite{Fabritiis2000}. In this approach, the SPH discretization of the Navier-Stokes equations is performed, and then thermal fluctuations are treated in the same way as in DPD. On the contrary,  the present meshfree method is a formulation which is based on a direct numerical discretization of the stochastic partial differential equation.
In the method the continuum constitutive model with  a stochastic stress tensor  is considered, and then a numerical approximation for the stochastic partial differential equations is employed.  An extension of the SDPD method including the  conservation of angular momentum has been presented by \citet{Mueller2015} to tackle  fluid problems where  angular momentum conservation is essential. Moreover, we note, that in  \cite{espanol2003} a  rotational friction force  governing particle spin interactions is included.
 In the present work, we  focus on the  Brownian motion of a particle due to  inherent fluctuations of the surrounding fluid.  Problems, where the conservation of angular momentum of the fluid particles is required, are left for future work. 
 The other important distinguishing feature of the present work, is the use of an incompressible fluid solver instead of a compressible one  as done in \cite{espanol2003, Huang2012, Poblete2014}. This allows treating the Brownian motion of a particle inside a liquid considered in the present work. 
 We note that a compressible fluid needs to be considered if one want to focus on the interaction between ultrasound waves and colloidal particles, as studied by \citet{Usabiaga2013}. For a  compressible solver for the fluctuating hydrodynamics equations in one dimension developed by the present authors we refer to  \cite{pandey2012}. The coupling of a suspended particle with a fluctuating compressible fluid is left to  future work.

In the present work we consider a fully Lagrangian meshfree particle method \cite{tiwari2002, tiwari2013}.  
The computational domain is approximated by moving grid points or particles. We note that a particle management procedure has to be
added in the method, see \cite{tiwari2002, tiwari2013} for details.
 The suitability of the method, for fluid-structure interaction with highly flexible structures in the case of regular flow fields has been shown by \citet{tiwari2007}. In this paper, we have extended this meshfree method to the coupling of rigid particles with  fluctuating fluids. 
 For validation, the Brownian motion of particles has been investigated. We have computed the velocity autocorrelation function (VACF) of the Brownian particle and compared it with the theoretical result, as given for example in \cite{Hinch1975}. A rigid sphere immersed in the incompressible fluctuating fluid has been considered to validate the numerical results.
\section{Governing Equations}
\label{GE}
We  consider a rigid sphere  inside an incompressible fluctuating fluid. 
Let $\Omega \subset \mathbb{R}^{3}$  denote the entire computational domain including both fluid and rigid body, the domain of the rigid body is denoted by $P$. A neutrally buoyant rigid particle is considered to demonstrate the Brownian motion of an immersed particle due to the inherent fluctuations in the fluid.

The governing equations for the motion of the incompressible fluctuating fluid are given by
\begin{equation}
\label{eqn17}
\frac{d\mathbf{x}}{dt} = \mathbf{u}  \qquad \text{in} \hspace{4 pt} \Omega\setminus P,
\end{equation}
\begin{equation}
\label{eqn1}
\rho_{f} \dfrac{d\mathbf{u}}{dt} = \nabla\cdot\boldsymbol{\sigma} \qquad\text{in} \hspace{4 pt} \Omega\setminus P,
\end{equation}
\begin{equation}
\label{eqn2}
\nabla \cdot \mathbf{u} = 0 \qquad \text{in} \hspace{4 pt} \Omega\setminus P,
\end{equation}
where 
$\mathbf{x}$ stands for the position vector of the fluid particle,
 $\rho_{f}$ denotes the density of the fluid. $\dfrac{d}{dt} = \dfrac{\partial}{\partial t} + \mathbf{u} . \nabla$ defines the material derivative. 
The stress tensor $\boldsymbol{\sigma}$  is given by
\begin{equation}
\label{eqn13}
\boldsymbol{\sigma} = -pI +\mu[\nabla\mathbf{u} + (\nabla\mathbf{u})^{T}]  + \mathbf{\widetilde{S}},
\end{equation}
where $p$ is the pressure and  $\mu$ is the dynamic viscosity of the surrounding fluid.
$\mathbf{\widetilde{S}}$ stands for the stochastic stress tensor, which models the inherent molecular fluctuations in the fluid.  The required stochastic properties  of $\mathbf{\widetilde{S}}$ have been derived  by Landau and Lifshitz \cite{Landau1959} in the spirit of a fluctuation-dissipation balance principle, described as
\begin{subequations}
\label{eqn14}
\begin{align}
\label{eqn14a}
&\langle \widetilde{S}_{ij}(\mathbf{x},t) \rangle = 0,\\
\label{eqn14b}
&\langle \widetilde{S}_{ik}(\mathbf{x},t) \widetilde{S}_{lm}(\mathbf{x^{'}},t^{'})\rangle = 2k_{B}T\mu(\delta_{il}\delta_{km} + \delta_{im}\delta_{kl})\delta(\mathbf{x} - \mathbf{x^{'}})\delta(t - t^{'}),
\end{align}
\end{subequations}
where $k_{B}$ is the Boltzmann constant, $T$ is the temperature of the fluid and $\langle ~\rangle$ is used for the ensemble averages. It has to be noted that originally these expressions have been derived for compressible fluids, but equation (\ref{eqn14}) is the corresponding approximation for an incompressible fluids. 

We note   that the non-linear LLNS equations define an ill-posed problem.  It has to be noticed that the stochastic forcing in the  LLNS equations is the divergence of a white noise process, rather than the more common external fluctuations modelled through white noise which have been discussed in \cite{Arnaud2009, Kloeden2009, Prato2004}. $\mathbf{\widetilde{S}}$ can not be defined pointwise either in space and time, therefore $\nabla\cdot\mathbf{\widetilde{S}}$ can not be given a precise mathematical interpretation. Further mathematical problems arise with  the interpretation of the non-linear term $\mathbf{u}\cdot\nabla\mathbf{u}$.
An approach to deal with these issues is to consider a regularization of the stochastic stress tensor, which is
typically the source of irregularity. The regularization can be physically justified by the fact that the fluctuating fields are defined from the underlying microscopic dynamics via spatial coarse-graining, as discussed in \cite{Espanol2009}. In such a   formalism, the nonlinear term and the stochastic forcing are naturally regularized by the discretization. 
We refer  to 
\citet{prato2002} for a   mathematical study of the two-dimensional Navier-Stokes equations  perturbed by a space-time white noise. 

The motion of the rigid sphere is governed by the Newton-Euler equations
\begin{equation}
\label{eqn6}
M \dfrac{d\mathbf{U}}{dt} = \mathbf{F}, 
\end{equation}
\begin{equation}
\label{eqn7}
\mathbf{I} \dfrac{d \boldsymbol\omega}{dt} = \mathbf{T} .
\end{equation}
Here, $\mathbf{U}$ and $\boldsymbol\omega$ represent the translational and rotational velocities of the sphere, respectively. $M$ and $\mathbf{I}$ denote mass and  moment of inertia of the rigid sphere, respectively. $\mathbf{F}$ is the resultant hydrodynamic force acting on the surface of the rigid sphere from the surrounding fluid,  
\begin{equation}
\label{eqn11}
\mathbf{F} = (-1)\int_{\partial P}\boldsymbol{\sigma} \textrm{\textbf{\^{n}}}_{s}\, \mathrm{d}(\partial P) .
\end{equation}
 $\mathbf{T}$ denotes the hydrodynamic torque of the hydrodynamical force 
\begin{equation}
\label{eqn12}
\mathbf{T} = (-1)\int_{\partial P}\mathbf{r}\times (\boldsymbol{\sigma}  \textrm{\textbf{\^{n}}}_{s})\, \mathrm{d}(\partial P),
\end{equation}
where $\textrm{\textbf{\^{n}}}_{s}$ is the unit outward normal on the surface of the sphere. $\mathbf{r} = \mathbf{x} - \mathbf{X}$ is the position vector with respect to the center of mass $(\mathbf{X})$ of the rigid body. The center of mass $\mathbf{X}$ and the orientation $\boldsymbol{\Theta}$ of the rigid body are updated by
\begin{equation}
\label{eqn8}
\dfrac{d\mathbf{X}}{dt} = \mathbf{U},
\end{equation}
\begin{equation}
\label{eqn9}
\dfrac{d\boldsymbol{\Theta}}{dt} = \boldsymbol\omega.
\end{equation}
Together, equations (\ref{eqn6} - \ref{eqn9}) describe the motion of the rigid body. No additional random term has been incorporated in the equation of motion of the immersed particle to model its Brownian motion. In this context, our approach is closer to \cite{sharma2006, patankar2004, Donev2014, Usabiaga2013}, but we have treated the fluid-structure system in a different manner.  

This formulation has to be complemented by appropriate initial and boundary conditions for the fluid-structure system. 
Let us denote the outer boundary of the computational domain $(\Omega)$, which is not shared by the rigid sphere, by $\Gamma$. 
We consider a simple cubic array of domains. Therefore periodic boundary conditions are employed on $\Gamma$, given as  
\begin{eqnarray}
\label{eqn3}
\nonumber &\mathbf{u}_{\Gamma_{L}} = \mathbf{u}_{\Gamma_{R}} \\ 
&\dfrac{\partial \mathbf{u}}{\partial \mathbf{\hat{n}}}\mid_{\Gamma_{L}} = - \dfrac{\partial \mathbf{u}}{\partial \mathbf{\hat{n}}}\mid_{\Gamma_{R}},
\end{eqnarray}
where $\mathbf{\hat{n}}$ denotes the unit outward normal on $\Gamma$. $\Gamma_{L}$ and $\Gamma_{R}$ are the left and right faces of the boundary $\Gamma$. Since we have periodicity in all directions, similar boundary conditions hold at the top - bottom and front - back faces. 

The no-slip boundary condition is considered on $\partial P$ which is the interior boundary of the fluid and the surface of the rigid body
\begin{equation}
\label{eqn4}
\mathbf{u} = \mathbf{U} + \boldsymbol\omega \times \mathbf{r} = \mathbf{v} \qquad\mathrm{on} \hspace{4 pt} \partial P .
\end{equation}
$\mathbf{v}$ represents the resultant velocity of the rigid sphere, which also gives the flow velocity at the interface.

The initial condition for the fluid is defined as 
\begin{equation}
\label{eqn5}
\mathbf{u}(t = 0) = \mathbf{u}_{0}, 
\end{equation}
where $\mathbf{u}_{0}$ should satisfy equation (\ref{eqn2}). The initial conditions for the suspended particle are
\begin{equation}
\label{eqn10}
\mathbf{X}(t = 0) = \mathbf{X}_{0}; \quad \boldsymbol{\Theta}(t = 0) = \boldsymbol{\Theta}_{0}; \quad \mathbf{U}(t = 0) = \mathbf{U}_{0}; \quad \boldsymbol\omega(t = 0) = \boldsymbol\omega_{0},
\end{equation} 
where $\mathbf{X}_{0}, \mathbf{U}_{0}, \boldsymbol\omega_{0}$ should satisfy equation (\ref{eqn4}) such that the resultant velocity of the rigid body $\mathbf{v}_{0}$ should be equal to the initial velocity $\mathbf{u}_{0}$ of the fluid.
Due to the consideration of a spherical particle made of homogeneous material, the non-linear term $\boldsymbol\omega \times \mathbf{I}\boldsymbol\omega$ vanishes in the equation (\ref{eqn7}). 

\section{Numerical Approximation}
In this section, we will discuss the numerical approximation of the coupled system.  
 
Though the occurrence of the divergence operator in front of the white noise makes the numerical approximation of the LLNS equations difficult, a systematic analysis of the numerical discretization of the LLNS equations in the context of finite-volume methods has been discussed in \cite{garcia2010}. We have also successfully simulated the compressible LLNS equation in one dimension in  a  meshfree framework \cite{pandey2012}. In the present work, this meshfree method is used for the spatial discretization of the incompressible LLNS equations. In this formulation, the discretization  in the spatial domain is defined by material points moving with the fluid velocity. These mesh particles carry all relevant physical properties of the fluid. 

For the time integration a projection based scheme for the motion of the fluid and an explicit Euler scheme for the motion of the rigid body have been used. We take a fixed time step $\Delta t$ throughout the computation. 

The spatio-temporal averaging of stochastic forcing is performed in the same way as explained in \cite{garcia2010, bell2012}. The components of the stochastic stresses are generated as 
\begin{equation}
\label{new18Feb1}
\widetilde{\mathbf{S}} = \left(\widetilde{S}_{ij}\right) = 
\begin{cases} \sqrt{\dfrac{4k_{B}T\mu}{\Delta V \Delta t}} \widetilde{\boldsymbol{\Re}} & \quad \text{for} ~i = j \\ \\ \sqrt{\dfrac{2k_{B}T\mu}{\Delta V \Delta t}} \widetilde{\boldsymbol{\Re}} & \quad \text{for} ~i \neq j
\end{cases},
\end{equation}
where $\widetilde{\boldsymbol{\Re}} = \dfrac{\boldsymbol{\Re} + \boldsymbol{\Re}^{T}}{2}$. A realization of $\widetilde{\boldsymbol{\Re}}$ is sampled using  a  stream of independent,  standard normally distributed  random numbers at each time step. $\Delta V = \dfrac{V}{N},$ where $V$ is the volume of the fluid domain $\Omega\setminus P$ and $N$ is the number of spatial discretization points at a particular time step.

It has to be noted that we update the value of $\Delta V$ at each time step dependent on the actual number of spatial points in the domain. In the meshfree formulation $N$ does not remain fixed over time due to  removing  spatial points which are too close to each other and adding new points in the sparse domain.  Therefore, $\Delta V$ changes over time and has to be appropriately updated at each time step. We note that  this approximation of the stochastic tensor has been successfully implemented in the meshfree numerical approximation of the compressible LLNS equations in one dimension, see \cite{pandey2012}.
\subsection{Time Integration}
At first we explain the  time integration of  the fluid equations and then the coupling of fluid and  rigid body motion.

\subsubsection{Time integration of the incompressible LLNS equations}
For the fluctuating hydrodynamic equations, we have employed an extension of the Chorin projection  scheme \cite{chorin1968}. In the next subsection, the spatial discretization of the differential operators will be described. These discretized differential operators are denoted by letter symbols to distinguish them from the corresponding continuum operators. For example, $\mathbf{G}$ denotes the discretized gradient operator, $\mathbf{L}$ and  $\mathbf{D}$ stand for the discretized Laplacian and divergence operators, respectively.The operation of a divergence operator on a tensor field such as the stochastic stress tensor $\widetilde{\mathbf{S}}$ is understood component-wise on the $x$, $y$, and $z$ components. Similarly, the gradient and Laplacian operators act component-wise on a vector. The consistency of the discretized differential operator has been shown for the solution of the incompressible Navier-Stokes equation, see \cite{tiwari2002}.

A temporal discretization of the  LLNS system has to reproduce the statistical properties of the continuum fluctuations, which is an additional challenge for the temporal scheme. We would like to mention that the temporal integration of fluctuating hydrodynamics can be higher order accurate only in the weak sense and only for the linearized equations of fluctuating hydrodynamics \cite{garcia2010}. 

We index the time step by a superscript $n$, i.e. quantities evaluated at time $n \Delta t $ are denoted by the superscript $n$.

At the first step of the projection scheme, we compute the new positions of the particles from equation (\ref{eqn17}) 
\begin{equation}
\label{eqn18}
\mathbf{x}^{n+1} = \mathbf{x}^{n} + \Delta t \mathbf{u}^{n}.
\end{equation}
At the new positions, we first compute the intermediate velocity $\mathbf{u^{*}}$ from equation (\ref{eqn1}) by ignoring the pressure term
\begin{equation}
\label{eqn19}
\mathbf{u}^{*} = \mathbf{u}^{n} + \dfrac{\Delta t}{\rho_{f}} (\mu\mathbf{L} \mathbf{u}^{n} + \mathbf{D} \widetilde{\mathbf{S}}^{n}).
\end{equation}
Then, we correct the intermediate velocity $\mathbf{u^{*}}$ by its projection onto the space of solenoidal velocity fields
\begin{equation}
\label{eqn20}
\mathbf{u}^{n+1} = \mathbf{u}^{*} - \dfrac{\Delta t}{\rho_{f}} \mathbf{G} p^{n+1}.
\end{equation}
This requires the gradient of the pressure field at the $(n+1)^{th}$ time step. For this a pressure Poisson equation is solved which comes from the incompressibility condition. Since
\begin{equation}
\label{new14Feb4}
\mathbf{D}\mathbf{u}^{n+1} = 0,
\end{equation}
we get the pressure Poisson equation from equation (\ref{eqn20})
\begin{equation}
\label{eqn22}
\mathbf{L} p^{n+1} = \dfrac{\rho_{f}}{\Delta t}\mathbf{D} \mathbf{u}^{*}.
\end{equation}
Thus we use the Neumann boundary condition for the pressure on the solid wall $\partial P$ which is given by 
 \begin{equation}
\label{eqn24}
\left(\mathbf{G}p^{n+1} \cdot \hat{\mathbf{n}_s}\right)|_{\partial P} =  0.
\end{equation}

We use the no-slip boundary condition for $\mathbf{u}^{*}$ and $\mathbf{u}^{n+1}$ on $\partial P$.
In this scheme, we update the positions of the mesh particles only at the beginning of every time step. Then the intermediate velocity, the final divergence free velocity field and the pressure field are  computed at these new particle positions. The stochastic flux $\mathbf{\widetilde{S}}$ is updated at each time step as described earlier. The differential operators appearing in equations (\ref{eqn19} -\ref{eqn24}) are computed at every particle position from its surrounding clouds of points. 
\subsubsection{Time integration for the solid structure}
For the fluid-structure interaction we have to couple the  Newton-Euler equations, given by equations (\ref{eqn6} - \ref{eqn9}), with the  LLNS equations. We use an explicit method for the time discretization  of the Newton-Euler equations. The main steps are as follows
\begin{itemize}
\item[1.] Once we get the value of $\mathbf{u}^{n+1},$ $p^{n+1}$ and the value of  $\widetilde{\mathbf{S}}^{n+1}$ (based on the stream of pseudo-random numbers ), compute the  stress $\boldsymbol{\sigma}^{n+1}$ according to equation (\ref{eqn13})
\begin{equation}
\label{eqn25}
\boldsymbol{\sigma}^{n+1} = -p^{n+1}I +\mu[\mathbf{G}\mathbf{u}^{n+1} + (\mathbf{G}\mathbf{u}^{n+1})^{T}]  + \widetilde{\mathbf{S}}^{n+1}.
\end{equation} 

\item[2.] Compute the hydrodynamic force and torque by a numerical approximation of the surface integrals, given by equations (\ref{eqn11}) and (\ref{eqn12}). To simulate the hydrodynamic interactions between the fluctuating fluid and the solid body, the solid body is defined by a boundary surface. In the numerical computation, this boundary surface is constructed by many point-like particles. 
We note that the surface of the sphere is triangulated using boundary particles. 
The triangulation of the sphere is illustrated in figure \ref{spherediscre}.
\begin{figure}[!h]
\centering
\includegraphics[scale=0.34]{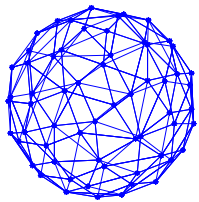}
\caption{Triangulation  to construct the boundary surface of a sphere.}
\label{spherediscre}       
\end{figure}
The components of the hydrodynamic force and the torque are computed as
\begin{equation}
\label{equ1_28Jan2014_1}
F_{i} = - \sum_{\mathbf{x} \in \partial P}\left( \left(\boldsymbol{\sigma}\cdot \mathbf{\hat{n}_{s}}\right)_{i}\right) _{\mathbf{x}}ds_{\mathbf{x}},
\end{equation}
\begin{equation}
\label{equ1_28Jan2014_2}
T_{i} = - \sum_{\mathbf{x} \in \partial P} \left(\left(\left(\mathbf{x - X}\right) \times \left(\boldsymbol{\sigma} \cdot  \mathbf{\hat{n}_{s}}\right) \right)_{i}\right)_{\mathbf{x}} ds_{\mathbf{x}},
\end{equation}

where $i$ runs from $1$ to $3,$ and $ds_{\mathbf{x}}$ is the area of the small surface element given by the triangulation.

\item[3.] Solve the equations of motion for the solid structure given by equations (\ref{eqn6}, \ref{eqn7} and \ref{eqn8}) together with  given initial conditions:
\begin{equation}
\label{eqn26}
\mathbf{X}^{n+1} = \mathbf{X}^{n} + \Delta t \mathbf{U}^{n},
\end{equation}
\begin{equation}
\label{eqn27}
\mathbf{U}^{n+1} = \mathbf{U}^{n} + \dfrac{\Delta t}{M}\mathbf{F}^{n},
\end{equation}
\begin{equation}
\label{eqn28}
\boldsymbol\omega^{n+1} = \boldsymbol\omega^{n} + \Delta t \mathbf{I^{-1}}T.
\end{equation}
\item[4.] Compute the resultant velocity of the structure which will yield the interface boundary condition for the fluid flow
\begin{equation}
\label{eqn29}
\mathbf{v}^{n+1}=\mathbf{u}^{n+1} = \mathbf{U}^{n+1} + \boldsymbol\omega^{n+1} \times  (\mathbf{x} - \mathbf{X}^{n+1}). 
\end{equation}
\item[5.]Move the rigid body with the resultant velocity $\mathbf{v}^{n+1}$.
Then, assign  new velocities to the interface boundary particles according to (\ref{eqn29}). The new positions of the boundary particles are obtained from the linear transformation representing the movement of the rigid sphere.
 \end{itemize}

It has been mentioned in  \cite{howard1992} that   fully explicit schemes for fluid-structure simulation, as the one used in the present implementation, will be unstable in certain situations due to the explicit discretization of the equations of translational motion of the rigid body.  
In \cite{howard1992} a condition for stability of the explicit scheme has been derived as $M_{v} < M$, where   $M_{v}$ denotes the so-called virtual mass of the fluid.

Now, in the case of a spherical object, the virtual mass in an infinite fluid medium is given by 
\begin{equation}
\label{eqn32}
M_{v} = \dfrac{2}{3} \pi r^{3} \rho_{f},
\end{equation} 
and the actual mass of the sphere is
\begin{equation}
\label{eqn33}
M = \dfrac{4}{3} \pi r^{3} \rho_{s}.
\end{equation} 
Hence, $M_{v} < M $ if $\rho_{f} < 2\rho_{s}$.\\
We have strictly followed this condition in our simulation by considering a neutrally buoyant spherical particle.

The discretizations described above yield a first order spatio-temporal discretization of the fluid-structure system. The time step $\Delta t$ is restricted by the CFL condition 
\begin{equation}
\label{eqn1_4May2015}
\Delta t \leq \textrm{min}\left\lbrace 0.16 \dfrac{h}{3U_{\textrm{max}}}, 0.11 \dfrac{h^{2}}{9\nu}\right\rbrace,
\end{equation}
where $\nu = \dfrac{\mu}{\rho_{f}}$ is the kinematic viscosity, $U_{max}$ is the typical advection speed and $h$ is the radius of the neighbourhood, which bounds the number of neighbour particles for a particular particle.  The time step is chosen such that none of the grid particles moves a distance more than the average spacing between them in a single time step. 
We note that we have chosen $h$ to be  equal to three times the initial spacing of the particles. Thus, the first term in
(\ref{eqn1_4May2015}) can be viewed as the usual hyperbolic CFL condition with a numerical factor  taking into account the 
possible changes of the grid spacing during the computation.
Similarly, the second term in (\ref{eqn1_4May2015}) is
the analogue to the parabolic CFL condition. Compare
\cite{Donev2014} for the corresponding conditions for a fixed grid.

This completes the temporal discretization of the fluid-structure system. The remaining task is to approximate the spatial derivatives.
\subsection{Spatial discretization}
An algorithm for the meshfree  solution of the incompressible Navier-Stokes equations  has been presented in \cite{tiwari2002}. 
Let $f  : \Omega \times [0,T] \longrightarrow \mathbb{R}$ be a scalar-valued function. $f(\mathbf{x}, t)$ denotes the value of the  function $f$ at position $\mathbf{x} \in \Omega$ at an instant $t$. In this formulation, the spatial differential operators acting on  $f$ are approximated at $\mathbf{x}$ in terms of the values of the function $f(\mathbf{x},t)$ at a set of neighbouring points of $\mathbf{x}$.  
We assign a weight function to each particle. This weight function is a function of the distance between the central particle (the particle on which the derivative is being computed) and its neighbouring particles. 
In order to restrict the number of neighbouring particles, we consider a weight function $w = w( \mathbf{x}_{j} - \mathbf{x}, h)$ with  compact support of size $h$. The choice of the weight function can be arbitrary. We choose a Gaussian weight function
\begin{equation}
\label{eqn34}
w\left(\mathbf{x}_{j} - \mathbf{x} ; h\right) =  \left\{ \begin{array}{ll} \exp\left( -\alpha\dfrac{\Vert \mathbf{x}_{j} - \mathbf{x} \Vert^{2}}{h^{2}}\right),  & \textrm{if $\dfrac{\Vert \mathbf{x}_{j} - \mathbf{x} \Vert}{h} \leq 1$} \\ 0, & \textrm{else} .
\end{array}
\right.
\end{equation} 
where $\alpha$ is a positive constant and $h$ defines the neighbourhood radius for $\mathbf{x}$.
 There are obvious restrictions on the distribution and the number of neighbouring particles of $\mathbf{x}$ in order to obtain a reasonable approximation of the derivatives.
  We remove  and add particles if the distribution is too dense or too sparse.
   The fluid quantities of newly added particles are approximated from their neighbouring values. Let $m$ be the number of neighbouring particles for the particle $\mathbf{x}$. Then, to approximate the spatial derivative of $f(\mathbf{x}, t)$ at $\mathbf{x} = (x, y, z)$ we consider the Taylor expansion of $\mathbf{f}(\mathbf{x}_{i}, t)$ around $\mathbf{x} = (x, y, z)$, for all neighbouring points of $\mathbf{x}$. This yields a linear system in the unknown derivatives, reading
\begin{equation}
\label{eqn36_29Jan2014}
\mathbf{e} = M\mathbf{a}-\mathbf{b}.
\end{equation}
where
\begin{equation}
\label{nonum1_29Jan2014}
M = \nonumber
\end{equation} 
\begin{equation}
{
\scriptstyle
\begin{bmatrix}
 \Delta x_{1} & \Delta y_{1} &\Delta z_{1} & \dfrac{\Delta x_{1}^{2}} {2} & \Delta x_{1}\Delta y_{1} & \Delta x_{1}\Delta z_{1}& \dfrac{\Delta y_{1}^{2}} {2}& \Delta y_{1}\Delta z_{1}&  \dfrac{\Delta z_{1}^{2}} {2} \\ 

\Delta x_{2} & \Delta y_{2} &\Delta z_{2} & \dfrac{\Delta x_{2}^{2}} {2} & \Delta x_{2}\Delta y_{2} & \Delta x_{2}\Delta z_{2}& \dfrac{\Delta y_{2}^{2}} {2}& \Delta y_{2}\Delta z_{2}&  \dfrac{\Delta z_{2}^{2}} {2} \\

\vdots & \vdots & \vdots & \vdots & \vdots & \vdots & \vdots & \vdots & \vdots \\

\Delta x_{m} & \Delta y_{m} &\Delta z_{m} & \dfrac{\Delta x_{m}^{2}} {2} & \Delta x_{m}\Delta y_{m} & \Delta x_{m}\Delta z_{m}& \dfrac{\Delta y_{m}^{2}} {2}& \Delta y_{m}\Delta z_{m}&  \dfrac{\Delta z_{m}^{2}} {2}\nonumber
\end{bmatrix},}
\end{equation}
\begin{align}
\label{nonum2_29Jan2014}
\mathbf{a} = 
\begin{array}{ccccccccc}
[f^{x} & f^{y} & f^{z} & f^{xx} & f^{xy} & f^{xz} & f^{yy} & f^{yz} & f^{zz}]^{T}, \nonumber
\end{array}
\\[7 pt]
\mathbf{b} = 
\begin{array}{ccccccccc}
[f_{1} -f & f_{2} - f & \hdots & \hdots & \hdots  & f_{m} - f]^{T}, \nonumber
\end{array}
\\[7 pt]
\mathbf{e} = 
\begin{array}{ccccccccc}
[e_{1} & e_{2} & \hdots & \hdots & \hdots  & e_{m}]^{T}. \nonumber
\end{array}
\end{align}
Here, $\Delta x_{i} = x_{i} - x, ~ \Delta y_{i} = y_{i} - y, ~ \Delta z_{i} = z_{i} - z$ for $i = 1, \hdots m$, and superscripts $x,~ y,~ z$ on $f $ represent the respective partial derivatives of the function.

This system is solved using the least square approximation by minimizing,
\begin{equation}
\label{eqn37_29Jan2014}
\mathit{J} =  \sum_{i = 1}^{m}w_{i}e_{i}^{2} = \left( M\mathbf{a} - \mathbf{b}\right)^{T} W \left(M\mathbf{a} - \mathbf{b}\right),
\end{equation} 
where, \\
$W = diag[w_{1}, \ldots, w_{m}]$ is a diagonal matrix, with entries $w_{i} = w\left(\mathbf{x}_{i} - \mathbf{x} ; h\right)$ as given in equation (\ref{eqn34}).

This yields the unknown $\mathbf{a}$ as
\begin{equation}
\label{eqn38_29Jan2014}
\mathbf{a} = \left(M^{T} W M\right)^{-1}\left(M^{T}W\right)\mathbf{b}.
\end{equation}
As a result, we will get derivatives of the prescribed function at a specific point as a linear combination of function values at its neighbour points.

In the present work, this approach will be also employed for the solution of the pressure Poisson equation with Neumann boundary conditions and for periodic boundary conditions of the velocity. It was first presented in \cite{tiwari2004} to solve the Poisson equation, where second-order convergence was demonstrated. The accuracy and stability of this method to solve the pressure Poisson equation with different boundary conditions was also discussed in \cite{tiwari2004}.

\section{Numerical results}
In this section we  present our simulation results to validate the numerical approximation, discussed in the previous sections.
We investigate the Brownian motion of a solid sphere due to the surrounding fluctuating fluid. In this context, the VACF of the sphere has been calculated and compared with the theoretical result as given in \cite{Hinch1975}.

To proceed, the governing equations are non-dimensionalized. The fundamental scales for the non-dimensionalization are chosen as
\begin{equation}
\label{eqn49}
\left.
\begin{array}{l l}
\text{Characteristic length} \longrightarrow \widetilde{\mathit{x}} = \Delta x \approx \dfrac{h}{3} \\\text{(average distance between the Lagrangian particles)},\\[10 pt]
\text{Characteristic time} \longrightarrow \widetilde{\mathit{t}} = \dfrac{\rho_{f} (\Delta x)^{2}}{\mu},\\[10 pt]
\text{Characteristic mass} \longrightarrow \widetilde{\mathit{M}} = \rho_{f} (\Delta x)^{3},
\end{array}\right\rbrace
\end{equation}

The resulting scalings for velocity, pressure and stress are
\begin{equation}
\label{eqn50}
\left.
\begin{array}{l l}
\text{Characteristic velocity} \longrightarrow \mathbf{\widetilde{\mathit{v}}} =  \dfrac{\widetilde{\mathit{x}}}{ \widetilde{\mathit{t}}} = \dfrac{\mu}{\rho_{f}\Delta x},\\[10 pt]

\text{Characteristic pressure, stress} \longrightarrow \widetilde{\mathit{P}}, \widetilde{\mathit{S}} = \dfrac{\mu^{2}}{\rho_{f}(\Delta x)^{2}}.
\end{array}\right\rbrace
\end{equation}

The resulting non-dimensionalized system is
\begin{equation}
\label{eqn1_27May2015}
\dfrac{d\mathbf{x}^{*}}{d t^{*}} =  \mathbf{u}^{*} \qquad\text{in} \quad \Omega\setminus P,
\end{equation}
\begin{equation}
\label{eqn51}
\dfrac{d\mathbf{u}^{*}}{d t^{*}} = -\nabla p^{*} + \triangle \mathbf{u^{*}} + \nabla \cdot \widetilde{\mathbf{S}}^{*} \qquad\text{in} \quad \Omega\setminus P,
\end{equation}
\begin{equation}
\label{eqn52}
\nabla \cdot \mathbf{u}^{*} = 0 \qquad \text{in} \quad \Omega\setminus P,
\end{equation}
where, properties of the random stresses $\widetilde{\mathbf{S}}^{*}$ are given as 
\begin{subequations}
\label{eqn14_1}
\begin{align}
\label{eqn14a_1}
&\langle \widetilde{S}_{ij}^{*}(\mathbf{x},t) \rangle = 0,\\
\label{eqn14b_1}
&\langle \widetilde{S}_{ik}^{*}(\mathbf{x},t) \widetilde{S}_{lm}^{*}(\mathbf{x^{'}},t^{'})\rangle = \dfrac{2k_{B}T\rho_{f}^{2}\widetilde{\mathit{x}}^{4}}{\mu^{3}}(\delta_{il}\delta_{km} + \delta_{im}\delta_{kl})\delta(\mathbf{x} - \mathbf{x^{'}})\delta(t - t^{'}),
\end{align}
\end{subequations}
\begin{equation}
\label{eqn54}
\dfrac{d\mathbf{U^{*}}}{dt^{*}} = \mathbf{F}^{*}\qquad \text{in} \quad P, 
\end{equation}
\begin{equation}
\label{eqn55}
\mathbf{I^{*}} \dfrac{d \boldsymbol\omega^{*}}{dt^{*}} = T^{*} \qquad \text{in} \quad P,
\end{equation}
\begin{equation}
\label{eqn56}
\dfrac{d\mathbf{X^{*}}}{dt^{*}} = \mathbf{U^{*}},
\end{equation}
\begin{equation}
\label{eqn57}
\dfrac{d\boldsymbol{\Theta}^{*}}{dt^{*}} = \boldsymbol\omega^{*},
\end{equation}
Non-dimensionalized variables are  denoted by the superscript $*$. 

\subsection{Brownian motion of a sphere in a three dimensional fluctuating fluid}
We consider a neutrally-buoyant spherical object immersed in an incompressible fluid. Initially, the sphere is placed at the center of a cubic domain. We use  periodic boundary conditions in all directions at the outer boundary $\Gamma$ and no-slip boundary conditions at the interface boundary $\partial P$. No additional force is applied other than the stochastic force. 
The only forces responsible for the motion of the sphere are the random stresses. 
Initially, the fluid-structure system is at rest. Now, we solve the system of equations (\ref{eqn1_27May2015}) - (\ref{eqn57}) for the three-dimensional fluid-structure system. We employ our scheme to perform a well known test \cite{patankar2004, atzberger2006, Iwashita2009, Donev2013, Donev2014} for the coupling of a spherical particle with the fluctuating fluid. In this context, we compute the VACF,
\begin{equation}
\label{eqn_11May2015}
C(t^{*}) = \dfrac{1}{d} \langle \mathbf{U}^{*}(t^{*}) \cdot \mathbf{U}^{*}(0)\rangle
\end{equation}
of a single Brownian particle diffusing through a periodic fluctuating fluid and compare the result with the analytical expression \cite{Hinch1975}. Here $d$ denotes the dimension of the computational domain. The VACF yields the crucial information about the Brownian dynamics at both short and long times. 

\begin{figure}[!h]
\centering
\includegraphics[scale=0.56]{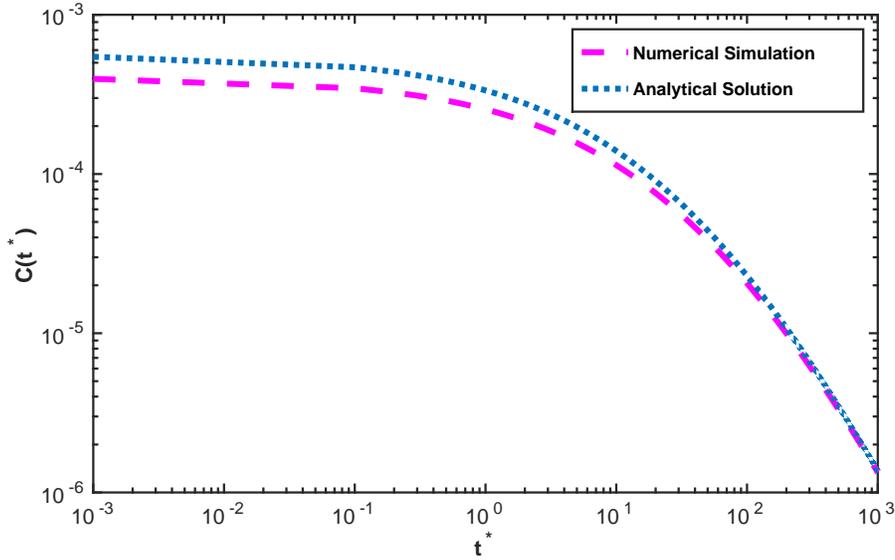}
\caption{The velocity autocorrelation function of a rigid sphere immersed in an incompressible fluctuating fluid.}
\label{fig_11May2015}       
\end{figure}

In figure \ref{fig_11May2015}, we show the VACF of a neutrally-buoyant   
spherical particle along with the analytical approximation of the VACF. We fill the cubic box with $N = 132350$  Lagrangian particles. The volume fraction of the spherical particle is \[\psi  = \dfrac{\text{volume of the spherical particle}}{\text{volume of the simulation box}} = 0.008.\] Table (\ref{tab1_12May2015}) gives the physical parameters used for this simulation. 
\begin{table}[!h]
\centering
\begin{tabular}{| c | c |}
\hline  
Fluid density $\rho_{f}$ & 1.0 \\
Viscosity $\mu$ & 1.0 \\
Thermal energy $k_{B}T$ & 0.83 \\
Hydrodynamic radius $R_{H}$ & 6.2\\
Schmidt number $S_{c}$ & 140.6 \\
\hline
\end{tabular}
\caption{Physical parameters.}
\label{tab1_12May2015}
\end{table}
All of these parameters are dimensionless, since they have been expressed via the fundamental scales for length, time etc., defined in equations (\ref{eqn49}) and (\ref{eqn50}). $R_H$ is the radius of the considered spherical particle. The Schmidt number $S_c$ is defined by the Stokes-Einstein diffusion coefficient $D$,
\[S_c \approx \dfrac{\nu}{D} = \dfrac{6 \pi \mu^{2} R_H}{\rho k_B T}\]

It can be observed that our numerical approximation correctly reproduces the long time behaviour of the Brownian dynamics. At the viscous time scale, $t^{*} > \tau_{\nu} = \rho_{f}R_{H}^{2}/\mu \approx 36$, the VACF shows the well-known algebraic decay $t^{-3/2}$ \cite{alder97}. There are some discrepancies in the result for short time dynamics. In particular, $C(t^{*})$ does not agree well with the analytical expression for $t^* < 10$. 
Such discrepancies in the short-time behaviour of the VACF have also been observed by other particle methods, as reported in \cite{Iwashita2009}. We note, however that a  very good approximation of the VACF, for the short and long-time behaviour, has been attained by \cite{Donev2014, Donev2013}, for a large range of Schmidt numbers using a Finite Volume approach.

In order to test further the validity of our method for the long-time behaviour of Brownian dynamics, we evaluate the long-time diffusion coefficient of the spherical particle. The long-time diffusion coefficient of the immersed particle is obtained by the discrete integral of the VACF,
\begin{equation}
\label{eqn1_12May2015}
D^{*}(t^{*}) = D^{*}(k \Delta t^{*}) \approx \dfrac{\Delta t^{*}}{2} C_{0} + \Delta t^{*} \sum_{i = 1}^{k-1}C_{i},
\end{equation} 
where $C_i$ are the values of the VACF
at time $t^{*} = i\Delta t^{*},$ for $i = 1 \ldots k-1$. 
We estimate the asymptotic value of the diffusion coefficient $D^{*} \approx D^{*}(t^{*}),$ for $t^{*} = L^{2}/ \nu,$ where $L$ denotes the size of the computational domain and $\nu$ is the dynamic viscosity of the fluid. The Stokes-Einstein relation for a periodic array of spheres in a cube can be written as 
\begin{equation}
\label{eqn2_12May2015}
D = \dfrac{k_{B}T}{6 \pi \mu R_{H} \xi_{c}},
\end{equation}
where $\xi_{c}$ denotes the correction factor to the drag coefficient of the immersed particle. In terms of non-dimensional variable, the drag coefficient correction factor can be re-written as 
\begin{equation}
\label{eqn2_15June2015}
\xi_{c} = \dfrac{k_{B}T\rho_{f}}{6 \pi \mu^{2} R_{H} D^{*}}.
\end{equation}
We have computed the drag coefficient factor from equation (\ref{eqn2_15June2015}) and compared it with the analytical value for different number of Lagrangian particles, given in table \ref{tab1}. The results shown in table \ref{tab1} are for the same volume fraction of the solid, which is $\psi  = 0.008$.
\begin{table}[!h]
\centering
\begin{tabular}{| c | c | c | c |}
\hline 
$M_a$ & Analytical  & Simulation Result & \% Error\\ & Value & Interval & \\\hline 
80000 & 1.525 & 1.4587  & 4.5 \\
100700 & 1.525 & 1.4899  & 2.3 \\
132350 & 1.525 & 1.5021  &1.5 \\
\hline
\end{tabular}
\caption{Comparison of numerical and analytical values of the drag correction factor.}
\label{tab1}
\end{table}
From table \ref{tab1} one can observe a  good agreement between the analytical and the numerical values of the drag correction factor. This confirms that the numerics is able to reproduce the Stokes-Einstein relation and to accurately represent the long-time Brownian dynamics of a solid sphere immersed in an incompressible fluid by solving the fluctuating hydrodynamic equations. 

\section{Conclusions}
A meshfree discretization for a system comprising a fluctuating incompressible fluid and a suspended solid body has been presented and validated via an investigation of the Brownian motion of the solid. No external forces on the solid structure other than the hydrodynamic forces from the surrounding fluid have been considered. The  LLNS equations for the fluctuating fluid have been coupled with the  Newton-Euler equations for the motion of the solid object.
 
Specifically, simulations of the Brownian motion of a sphere have been performed. 
To validate the simulation, the VACF of the Brownian particle has been calculated and compared with the available analytical expression. The long-time Brownian dynamics of the sphere has also been confirmed by the Stokes-Einstein relation. In that context, we have calculated the corresponding correction factor for the drag coefficient in a three-dimensional periodic system. 
The numerical value of the correction factor shows good agreement with the analytical result. 

The fact that the Stokes-Einstein relation could be reproduced demonstrates that the method is consistent with the continuum fluctuation-dissi-\\pation theorem (FDT). 
This still needs to be verified at the discrete level. The numerical validation of the discrete FDT by static and dynamic structure factors is a future task of the authors. 
Moreover, a detailed study of the problem for a larger range of Schmidt numbers, as given in \citep{Donev2013}, has  also to be considered. 
The study of fluid-fluid interfaces and the dynamics of small particles at these interfaces with fluctuating hydrodynamics in the framework of a meshfree discretization is another subject of  future work.

\section*{Acknowledgment:}  This work was partially supported by the German Research Foundation (DFG), grant number KL 1105/20-1, RTG 1932 ``Stoch-astic models in the engineering sciences''  and by the German Academic Exchange Service (DAAD), PhD program MIC. 

\bibliographystyle{model1-num-names}
\bibliography{myref}

\end{document}